\documentclass[twocolumn]{IEEEtran}
\usepackage[cmex10]{amsmath}
\usepackage{amsmath}
\usepackage{amssymb}

\usepackage{graphicx}          
\usepackage[dvips]{epsfig}    
\newtheorem{thm}{\bf Theorem}[section]
\newtheorem{assume}[thm]{\bf Assumption}
\newtheorem{problem}[thm]{\bf Problem}
\newtheorem{remark}[thm]{\bf Remark}
\newtheorem{lemma}[thm]{\bf Lemma}
\newtheorem{proposition}[thm]{\bf Proposition}
\newtheorem{definition}[thm]{\bf Definition}
\usepackage{tikz}
\usetikzlibrary{shapes,arrows}

\title{Global Adaptive Dynamic Programming for Continuous-Time Nonlinear Systems}
\author{Yu~Jiang,~\IEEEmembership{Member,~IEEE}, Zhong-Ping~Jiang,~\IEEEmembership{Fellow,~IEEE}
\thanks{This work has been supported in part by the National Science Foundation under Grants ECCS-1101401 and ECCS-1230040. The work was done when Y. Jiang was a PhD student in the Control and Networks Lab at New York University Polytechnic School of Engineering.}
\thanks{Y. Jiang is with The MathWorks, 3 Apple Hill Dr, Natick, MA 01760, E-mail: yu.jiang@mathworks.com} 
\thanks{Z. P. Jiang is with the Department
of Electrical and Computer Engineering, Polytechnic School of Engineering, New York University, Five Metrotech Center, Brooklyn,
NY 11201 USA. zjiang@nyu.edu}
}

\begin{document}

%


%

\maketitle

\begin{abstract}                          
This paper presents a novel method of global adaptive dynamic programming (ADP) for the adaptive optimal control of nonlinear polynomial systems. The  strategy consists of relaxing the problem of solving the Hamilton-Jacobi-Bellman (HJB) equation to an optimization problem, which is solved via a new policy iteration method. The proposed method distinguishes from previously known nonlinear ADP methods in that the neural network approximation is avoided, giving rise to significant computational improvement. Instead of semiglobally or locally stabilizing, the resultant control policy is globally stabilizing for a general class of nonlinear polynomial systems.  Furthermore, in the absence of the a priori knowledge of the system dynamics, an online learning method is devised to implement the proposed policy iteration technique by generalizing the current ADP theory. Finally, three numerical examples are provided to validate the effectiveness of the proposed method.
\end{abstract}

\begin{IEEEkeywords}
Adaptive dynamic programming, nonlinear systems, optimal control, global stabilization.
\end{IEEEkeywords}


\section{Introduction}
Dynamic programming \cite{bellman1957dynamic} offers a theoretical way to solve optimal control problems. However, it suffers from the inherent computational complexity, also known as the {\it curse of dimensionality}. Therefore, the need for approximative methods has been recognized as early as in the late 1950s \cite{bellman1959functional}. Within all these methods, adaptive/approximate dynamic programming (ADP) \cite{bertsekas2007dynamic, bertsekas1996neuro, powell2007approximate, si2004handbook, werbos1974beyond, lewisliubook2013} is a class of heuristic techniques that solves for the cost function by searching for a suitable approximation. In particular, adaptive dynamic programming \cite{werbos1974beyond, werbos1977advanced} employs the idea from reinforcement learning \cite{sutton1998reinforcement} to achieve online approximation of the cost function, without using the knowledge of the system dynamics.  ADP has been extensively studied for Markov decision processes (see, for example, \cite{bertsekas1996neuro} and  \cite{powell2007approximate}), as well as dynamic systems (see the review papers \cite{lewis2009reinforcement} and \cite{wang2009adaptive}). Stability issues in ADP-based control systems design are addressed in \cite{balakrishnan2008issues, lewisliubook2013fullbook, vrabie2013book}. A robustification of ADP, known as Robust-ADP or RADP, is recently developed by taking into account dynamic uncertainties \cite{jiang2013robust}.

To achieve online approximation of the cost function and the control policy, neural networks are widely used in the previous ADP architecture. Although neural networks can be used as universal approximators \cite{hornik1989multilayer}, \cite{park1991universal}, there are at least two major limitations for ADP-based online implementations. First, in order to approximate unknown functions with high accuracy, a large number of basis functions comprising the neural network are usually required. Hence, it may incur a huge computational burden for the learning system. Besides, it is not trivial to specify the type of basis functions, when the target function to be approximated is unknown. Second, neural network approximations generally are effective only on some compact sets, but not in the entire state space. Therefore, the resultant control policy may not provide global asymptotic stability for the closed-loop system. In addition, the compact set, on which the uncertain functions of interest are to be approximated, has to be carefully quantified before one applies the online learning method, such that stability can be assured during the learning process \cite{jiang2012tnnls}.

The main purpose of this paper is to develop a novel ADP methodology to achieve global and adaptive suboptimal stabilization of uncertain continuous-time nonlinear polynomial systems via online learning. As the first contribution of this paper, an optimization problem, of which the solutions can be easily parameterized, is proposed to relax the problem of solving the Hamilton-Jacobi-Bellman (HJB) equation. This approach is inspired from the relaxation method used in approximate dynamic programming for Markov decision processes (MDPs) with finite state space \cite{de2003linear}, and more generalized discrete-time systems \cite{lincoln2006relaxing, wang2010approximate, wang2011performance, savorgnan2009discrete, schweitzer1985generalized, summers2012approximate}. However, methods developed in these papers cannot be trivially extended to the continuous-time setting, or achieve global asymptotic stability of general nonlinear systems. The idea of relaxation was also used in nonlinear $\mathcal{H}_{\infty}$ control, where Hamilton-Jacobi inequalities are used for nonadaptive systems \cite{helton1999extending, schaft1999l2, Sassano2012dynamic}.

The second contribution of the paper is to propose a relaxed policy iteration method. Each iteration step is formulated as a sum of squares (SOS) program \cite{parrilo2000structured, blekherman2013semidefinite}, which is equivalent to a semidefinite program (SDP). It is worth pointing out that, different from the inverse optimal control design \cite{krstic1998inverse},  the proposed method finds directly a suboptimal solution to the original optimal control problem.

The third contribution is an online learning method that implements the proposed iterative schemes using only the real-time online measurements, when the perfect system knowledge is not available. This method can be regarded as a nonlinear variant of our recent work for continuous-time linear systems with completely unknown system dynamics \cite{jiang2012computational}. This method distinguishes from previously known nonlinear ADP methods in that the neural network approximation is avoided for computational benefits and that the resultant suboptimal control policy achieves global asymptotic stabilization, instead of only semi-global or local stabilization.

The remainder of this paper is organized as follows. Section 2 formulates the problem and introduces some basic results regarding nonlinear optimal control. Section 3 relaxes the problem of solving an HJB equation to an optimization problem. Section 4 develops a relaxed policy iteration technique for polynomial systems using sum of squares (SOS) programs \cite{blekherman2013semidefinite}. Section 5 develops an online learning method for applying the proposed policy iteration, when the system dynamics are not known exactly.  Section 6 examines three numerical examples to validate the efficiency and effectiveness of the proposed method. Section 7 gives concluding remarks.

{\it Notations:} Throughout this paper, we use $\mathcal{C}^1$ to denote the set of all continuously differentiable functions. $\mathcal{P}$ denotes the set of all functions in $\mathcal{C}^1$ that are also positive definite and proper. $\mathbb{R}_+$ indicates the set of all non-negative real numbers. For any vector $u\in\mathbb{R}^m$ and any symmetric matrix $R\in\mathbb{R}^{m\times m}$, we define $|u|_R^2$ as $u^TRu$.
A feedback control policy $u$ is called globally stabilizing, if under this control policy, the closed-loop system is globally asymptotically stable  at the origin \cite{khalil2002nonlinear}.
Given a vector of polynomials $f(x)$, ${\rm deg}(f)$ denotes the highest polynomial degree of all the entries in $f(x)$.
For any positive integers $d_1, d_2$ satisfying $d_2\ge d_1$, $\vec{m}_{d_1,d_2}(x)$ is the vector of all $(_{~d_2}^{n+d_2})-(_{~d_1-1}^{n+d_1-1})$ distinct monic monomials in $x\in\mathbb{R}^n$
with degree at least $d_1$ and at most $d_2$, and arranged in lexicographic order \cite{cox2007ideals}. Also, $\mathbb{R}[x]_{d_1,d_2}$ denotes the set of all polynomials in $x\in\mathbb{R}^n$ with degree no less than $d_1$ and no greater than $d_2$.
$\nabla V$ refers to the gradient of a function $V:\mathbb{R}^n\rightarrow\mathbb{R}$.

\section{Problem formulation and preliminaries}
\label{sec:Prob_form}
\subsection{Problem formulation}
Consider the nonlinear system
\begin{eqnarray}
\dot{x}=f(x)+g(x)u\label{eq:sys1}
\end{eqnarray}
where $x\in\mathbb{R}^{n}$ is the system state, $u\in\mathbb{R}^{m}$
is the control input, $f:\mathbb{R}^n\rightarrow\mathbb{R}^n$ and $g:\mathbb{R}^n\rightarrow\mathbb{R}^{n\times m}$ are polynomial mappings with $f(0)=0$.


In conventional optimal control theory \cite{lewis2012optimal}, the common objective is to find a control policy $u$ that
minimizes certain performance index. In this paper, it is specified as follows.
\begin{eqnarray}
J(x_0,u)=\int_{0}^{\infty}r(x(t),u(t))dt,~~x(0)=x_0 \label{eq:cost}
\end{eqnarray}
where
$r(x,u)=q(x)+u^{T}R(x)u$, with $q(x)$ a positive definite polynomial function, and $R(x)$ is a symmetric positive definite matrix of polynomials. Notice that, the purpose of specifying $r(x,u)$ in this form is to guarantee that an optimal control policy can be explicitly determined, if it exists.


\begin{assume}
\label{ass:controlability}
Consider system \eqref{eq:sys1}.
There exist a function $V_0\in\mathcal{P}$ and a feedback control policy $u_{1}$, such that
\begin{eqnarray}
\mathcal{L}(V_0(x),u_1(x))\ge 0,~~\forall x\in\mathbb{R}^n\label{eq:ini}
\end{eqnarray}
where, for any $V\in\mathcal{C}^1$ and $u\in\mathbb{R}^m$,
\begin{eqnarray}
\mathcal{L}(V,u)=-\nabla V^T(x)(f(x)+g(x)u)-r(x,u). \label{eq:def_L}
\end{eqnarray}
\end{assume}

Under Assumption \ref{ass:controlability}, the closed-loop system composed of \eqref{eq:sys1} and $u=u_1(x)$ is globally asymptotically stable at the origin, with a well-defined Lyapunov function $V_0$. With this property, $u_1$ is also known as an {\it admissible} control policy \cite{beard1997galerkin}, implying that the cost $J(x_0,u_1)$ is finite, $\forall x_0\in\mathbb{R}^n$. Indeed, integrating both sides of \eqref{eq:ini} along the trajectories of the closed-loop system composed of \eqref{eq:sys1} and $u=u_1(x)$ on the interval $[0,+\infty)$, it is easy to show that
\begin{eqnarray}
J(x_0,u_1)\le V_0(x_0),~~\forall x_0\in\mathbb{R}^n.
\end{eqnarray}

\subsection{Optimality and stability}
\label{sec:optnstb}
Here, we recall a basic result connecting optimality and global asymptotic stability in nonlinear systems \cite{SepulchreBook}.

To begin with, let us give the following assumption.

\begin{assume}
\label{ass:HJB}
There exists $V^{\rm o}\in\mathcal{P}$,
such that the Hamilton-Jacobi-Bellman (HJB) equation holds
\begin{eqnarray}
\mathcal{H}(V^{\rm o})=0 \label{eq:HJB}
\end{eqnarray}
where
\begin{eqnarray*}
\mathcal{H}(V)&=&\nabla V^T(x)f(x)+q(x) \\
&& -\frac{1}{4}\nabla V^T(x)g(x)R^{-1}(x)g^T(x)\nabla V(x).
\end{eqnarray*}
\end{assume}

Under Assumption \ref{ass:HJB}, it is easy to see that $V^{\rm o}$ is a well-defined Lyapunov function for the closed-loop system comprised of \eqref{eq:sys1} and
\begin{eqnarray}
u^{\rm o}(x)=-\frac{1}{2}R^{-1}(x)g^T(x)\nabla V^{\rm o}(x). \label{eq:optimality}
\end{eqnarray}

Hence, this closed-loop system is globally asymptotically stable at $x=0$ \cite{khalil2002nonlinear}. Then, according to \cite[Theorem 3.19]{SepulchreBook},
$u^{\rm o}$ is the optimal control policy, and the value function $V^{\rm o}(x_0)$ gives the optimal cost at the initial condition $x(0)=x_0$, i.e.,
\begin{eqnarray}
V^{\rm o}(x_0)=\min_u J(x_0,u)=J(x_0,u^{\rm o}),~~\forall x_0\in\mathbb{R}^n. \label{eq:optimal_cost}
\end{eqnarray}

It can also be shown that $V^{\rm o}$ is the unique solution to the HJB equation \eqref{eq:HJB} with $V^{\rm o}\in\mathcal{P}$. Indeed, let $\hat V\in\mathcal{P}$ be another solution to \eqref{eq:HJB}. Then, by Theorem 3.19 in \cite{SepulchreBook}, along the solutions of the closed-loop system composed of \eqref{eq:sys1} and $u=\hat u=-\frac{1}{2}R^{-1}g^T\nabla \hat V$, it follows that
\begin{eqnarray}
\hat V(x_0) = V^{\rm o}(x_0)-\int_0^{\infty}|u^{\rm o}-\hat u|_R^2dt, 
            ~~~~~~\forall x_0\in\mathbb{R}^n. \label{eq:optimality2}
\end{eqnarray}

Finally, by comparing \eqref{eq:optimal_cost} and \eqref{eq:optimality2}, we conclude that 
$$V^{\rm o}=\hat V.$$

\subsection{Conventional policy iteration}

The above-mentioned result implies that, if there exists a class-$\mathcal{P}$ function which solves the HJB equation \eqref{eq:HJB}, an optimal control policy can be obtained.
However, the nonlinear HJB equation \eqref{eq:HJB} is very difficult to be solved analytically in general. As a result, numerical methods are developed to approximate the solution. In particular, the following policy iteration method is widely used \cite{saridis1979approximation}.

\begin{enumerate}
\item[1)] {\it Policy evaluation:}~For $i=1,2,\cdots$, solve for the cost function $V_i(x)\in\mathcal{C}^1$, with $V_i(0)=0$, from the following partial differential equation.
\begin{eqnarray}
\mathcal{L}(V_i(x),u_i(x))=0.\label{eq:GHJB}
\end{eqnarray}

\item[2)] {\it Policy improvement:} Update the control policy by
\begin{eqnarray}
u_{i+1}(x)=-\frac{1}{2}R^{-1}(x)g^{T}(x)\nabla V_{i}(x). \label{eq:pimp}
\end{eqnarray}
\end{enumerate}

The convergence property of the conventional policy iteration algorithm is given in the following theorem, which is an extension of \cite[Theorem 4]{saridis1979approximation}. For the readers' convenience, its proof is given in Appendix \ref{sec:proofConv}.

\begin{thm}
\label{thm:conv}
Suppose Assumptions \ref{ass:controlability} and  \ref{ass:HJB} hold, and the solution $V_i(x)\in\mathcal{C}^1$ satisfying  \eqref{eq:GHJB} exists, for $i=1,2,\cdots$. Let $V_i(x)$ and $u_{i+1}(x)$ be the functions generated from (\ref{eq:GHJB}) and (\ref{eq:pimp}). Then,  the following properties hold, $\forall i=0,1,\cdots$.
\begin{enumerate}
\item[1)] $V^{\rm o}(x)\le V_{i+1}(x)\le V_i(x)$, $\forall x\in \mathbb{R}^n$;
\item[2)] $u_{i+1}$ is globally stabilizing;
\item[3)] suppose there exist $V^*\in\mathcal{C}^1$ and $u^*$, such that $\forall x_0\in\mathbb{R}^n$, we have $\lim\limits_{i\rightarrow\infty}V_i(x_0)=V^*(x_0)$ and $\lim\limits_{i\rightarrow\infty}u_i(x_0)=u^*(x_0)$. Then, $V^*=V^{\rm o}$ and $u^*=u^{\rm o}$.
\end{enumerate}
\end{thm}

Notice that finding the analytical solution to \eqref{eq:GHJB} is still non-trivial. Hence, in practice, the solution is approximated using, for example, neural networks or Galerkin's method \cite{beard1997galerkin}. Also, ADP-based online approximation method can be applied to compute numerically the cost functions via online data \cite{vrabie2009neural}, \cite{jiang2012tnnls}, when the precise knowledge of $f$ or $g$ is not available. However, although approximation methods can give acceptable results on some compact set in the state space, they cannot be used to achieve global stabilization. In addition, in order to reduce the approximation error, huge computational complexity is almost inevitable.

\section{Suboptimal control with relaxed HJB equation} 
%


In this section, we consider an auxiliary optimization problem, which allows us to obtain a suboptimal solution to the
minimization problem
\eqref{eq:cost}
subject to \eqref{eq:sys1}.
For simplicity, we will omit the arguments of functions whenever there is no confusion in the context.
\begin{problem}[Relaxed optimal control problem]
\label{prob:relaxed_HJB}
\begin{eqnarray}
\min_{V}&\int_\Omega V(x)dx \label{eq:wcost}\\
{\rm s.t.}~~&\mathcal{H}(V)\le 0&\label{eq:IE}\\
&V\in\mathcal{P}&
\end{eqnarray}
where $\Omega\subset\mathbb{R}^n$ is an arbitrary compact set containing the origin as an interior point. As a subset of the state space, $\Omega$ describes the area in which the system performance is expected to be improved the most.
\end{problem}

\begin{remark}
Notice that Problem \ref{prob:relaxed_HJB} is called a {\it relaxed} problem of \eqref{eq:HJB}. Indeed, if we restrict this problem by replacing the inequality constraint \eqref{eq:IE} with the equality constraint \eqref{eq:HJB}, there will be only one feasible solution left and the objective function can thus be neglected. As a result, Problem \ref{prob:relaxed_HJB} reduces to the problem of solving \eqref{eq:HJB}.
\end{remark}


Some useful facts about Problem \ref{prob:relaxed_HJB} are given as follows.
\begin{thm}
\label{thm:relaxed_HJB}
Under Assumptions \ref{ass:controlability} and \ref{ass:HJB}, the following hold.
 \begin{enumerate}
\item[1)] Problem \ref{prob:relaxed_HJB} has a nonempty feasible set.
\item[2)] Let $V$ be a feasible solution to Problem \ref{prob:relaxed_HJB}. Then, the control policy
\begin{eqnarray}
\bar u(x)=-\frac{1}{2}R^{-1}g^T\nabla V\label{eq:controller_v}
\end{eqnarray}
 is globally stabilizing.
\item[3)] For any $x_0\in\mathbb{R}^n$, an upper bound of the cost of the closed-loop system comprised of \eqref{eq:sys1} and \eqref{eq:controller_v} is given by $V(x_0)$, i.e.,
\begin{eqnarray}
       J(x_0,\bar u) \le V(x_0).
\end{eqnarray}
\item[4)] Along the trajectories of the closed-loop system \eqref{eq:sys1} and \eqref{eq:optimality}, the following inequalities hold for any $x_0\in\mathbb{R}^n$:
\begin{eqnarray}
V(x_0)+\int_0^{\infty}\mathcal{H}(V(x(t))) dt
\le V^{\rm o}(x_0)\le V(x_0). \label{eq:V_ieq}
\end{eqnarray}
\item[5)] $V^{\rm o}$ as defined in \eqref{eq:optimal_cost} is a global optimal solution to Problem \ref{prob:relaxed_HJB}.
\end{enumerate}
\end{thm}

\begin{IEEEproof}
%

1) Define  $\displaystyle u_0=-\frac{1}{2}R^{-1}g^T\nabla V_0$. Then,
\begin{eqnarray*}
\mathcal{H}(V_0)&=&\nabla V_0^T\left(f+gu_0\right)+r(x,u_0)\\
&=&\nabla V_0^T\left(f+gu_1\right)+r(x,u_1)\\
&&+\nabla V_0^Tg(u_0-u_1)+|u_0|^2_R-|u_1|^2_R\\
&=&\nabla V_0^T\left(f+gu_1\right)+r(x,u_1)\\
&&-2|u_0|^2_R-2u_0^TRu_1+|u_0|^2_R-|u_1|^2_R\\
&=&\nabla V_0^T\left(f+gu_1\right)+r(x,u_1)-|u_0-u_1|^2_R\\
&\le&0
\end{eqnarray*}

Hence, $V_0$ is a feasible solution to Problem \ref{prob:relaxed_HJB}.

2) To show global asymptotic stability, we only need to prove that $V$ is a well-defined Lyapunov function for the closed-loop system composed of \eqref{eq:sys1} and \eqref{eq:controller_v}. Indeed, along the solutions of the closed-loop system, it follows that
\begin{eqnarray*}
\dot V &=&\nabla V^T(f+g\bar u) \\
      &=&  \mathcal{H}(V)-r(x,\bar u) \\
      &\le& -q(x)
\end{eqnarray*}
Therefore, the system is globally asymptotically stable at the origin \cite{khalil2002nonlinear}.

3) Along the solutions of the closed-loop system comprised of \eqref{eq:sys1} and \eqref{eq:controller_v}, we have
\begin{eqnarray}
V(x_0)&=&-\int_0^{T}\nabla V^T(f+g\bar u)dt+V(x(T))\nonumber\\
&=&~\int_0^{T}\left[r(x,\bar u)-\mathcal{H}(V)\right]dt+V(x(T))\nonumber\\
&\ge&\int_0^{T}r(x,\bar u)dt+V(x(T))\label{eq:int_T}
\end{eqnarray}

By 2), $\lim\limits_{T\rightarrow+\infty}V(x(T))=0$. Therefore, letting $T\rightarrow+\infty$, by \eqref{eq:int_T} and \eqref{eq:cost}, we have
\begin{eqnarray}
 V(x_0)\ge J(x_0,\bar u).
 \end{eqnarray}

4) By 3), we have
\begin{eqnarray}
 V(x_0)\ge J(x_0,{\bar u})\ge \min\limits_{u}J(x_0,{\bar u})=V^{\rm o}(x_0).
 \end{eqnarray}
Hence, the second inequality in \eqref{eq:V_ieq} is proved.

On the other hand,
\begin{eqnarray*}
\mathcal{H}(V)&=&\mathcal{H}(V)-\mathcal{H}(V^{\rm o})\\
&=&(\nabla  V-\nabla V^{\rm o})^T(f+g{\bar u})+r(x,{\bar u})\\
&&- (\nabla V^{\rm o})^Tg(u^{\rm o}-{\bar u})-r(x,u^{\rm o})\\
&=&(\nabla V-\nabla V^{\rm o})^T(f+gu^{\rm o})-|{\bar u}-u^{\rm o}|^2_R\\
&\le&(\nabla V-\nabla V^{\rm o})^T(f+gu^{\rm o})
\end{eqnarray*}

Integrating the above equation along the solutions of the closed-loop system \eqref{eq:sys1} and \eqref{eq:optimality} on the interval $[0,+\infty)$, we have
\begin{eqnarray}
V(x_0)-V^{\rm o}(x_0)\le -\int_0^{\infty}\mathcal{H}(V(x))dt.
\end{eqnarray}

5) By 3), for any feasible solution $V$ to Problem \ref{prob:relaxed_HJB}, we have $V^{\rm o}(x)\le V(x)$. Hence,
\begin{eqnarray}
\int_{\Omega}V^{\rm o}(x)dx\le \int_{\Omega}V(x)dx
\end{eqnarray}
which implies that $V^{\rm o}$ is a global optimal solution.

The proof is therefore complete.
\end{IEEEproof}

\begin{remark}
A feasible solution $V$ to Problem \ref{prob:relaxed_HJB} may not necessarily be the true cost function associated with the control policy $\bar u$ defined in \eqref{eq:controller_v}. However, by Theorem \ref{thm:relaxed_HJB}, we see $V$ can be viewed as an upper bound or an {\it overestimate} of the actual cost, inspired by the concept of {\it underestimator} in \cite{wang2010approximate}. Further,  $V$ serves as a Lyapunov function for the closed-loop system and can be more easily parameterized than the actual cost function. For simplicity, $V$ is still called the {\it cost function}, in the remainder of the paper.
\end{remark}

\section{SOS-based Policy Iteration for polynomial systems}

The inequality constraint \eqref{eq:IE} contained in Problem \ref{prob:relaxed_HJB} provides us the freedom of specifying desired analytical forms of the cost function. However, solving \eqref{eq:IE} is non-trivial in general, even for polynomial systems (see, for example, \cite{chen2002global}, \cite{franze2012constrained}, \cite{moulay2005stabilization}, \cite{sontag1979observability}, \cite{xu2009simultaneous}). Indeed, for any polynomial with degree no less than four, deciding its non-negativity is an NP-hard problem \cite{parrilo2000structured}. Thanks to the developments in sum of squares (SOS) programs \cite{blekherman2013semidefinite, parrilo2000structured}, the computational burden can be significantly reduced, if inequality constraints can be restricted to SOS constraints. The purpose of this section is to develop a novel policy iteration method for polynomial systems using SOS-based methods \cite{blekherman2013semidefinite, parrilo2000structured}.

\subsection{Polynomial parametrization}

Let us first assume $R(x)$ is a constant, real symmetric matrix. Notice that $\mathcal{L}(V_i,u_i)$ is a polynomial in $x$,  if $u_i$ is a polynomial control policy and $V_i$ is a polynomial cost function. Then,  the following implication holds
\begin{eqnarray}
\mathcal{L}(V_i,u_i){~\rm~is~SOS} \Rightarrow \mathcal{L}(V_i,u_i)\ge 0.
\end{eqnarray}

In addition, for computational simplicity, we would like to find some positive integer $r$,  such that $V_i\in \mathbb{R}[x]_{2,2r}$. Then, the new control policy $u_{i+1}$ obtained from \eqref{eq:sosu} is a vector of polynomials.

Based on the above discussion, the following Assumption is given to replace Assumption \ref{ass:controlability}.
\begin{assume}
\label{ass:polyall}
There exist smooth mappings $V_0:\mathbb{R}^n\rightarrow
\mathbb{R}$ and $u_1:\mathbb{R}^n\rightarrow \mathbb{R}^{m}$, such that $V_0\in\mathbb{R}[x]_{2,2r}\cap\mathcal{P}$ and $\mathcal{L}(V_0,u_1)$ is SOS.
\end{assume}

%
%

\subsection{SOS-programming-based policy iteration}

Now, we are ready to propose a relaxed policy iteration scheme.
Similar as in other policy-iteration-based iterative schemes, an initial globally stabilizing (and admissible) control policy has been assumed in Assumption \ref{ass:polyall}.

\begin{enumerate}
\item[1)] {\it Policy evaluation:} For $i=1,2,\cdots$, solve for an optimal solution $p_i$ to the following optimization program:
\begin{eqnarray}
\min_{p} &&\int_{\Omega}V(x)dx
\label{eq:sosopt1}\\
{\rm s.t.}~~~~\mathcal{L}(V,u_i)&{\rm is}&{\rm SOS} \label{eq:sosopt3}\\
V_{i-1}-V&{\rm is}&{\rm SOS}    \label{eq:sosopt4}
%
%
 %
\end{eqnarray}
where $V=p^T\vec{m}_{2,2r}(x)$. Then, denote $V_i = p_i^T\vec{m}_{2,2r}(x)$.

\item[2)] {\it Policy improvement:} Update the control policy by
\begin{eqnarray}
u_{i+1}=-\frac{1}{2}R^{-1}g^{T}{ \nabla V_i}. \label{eq:sosu}
\end{eqnarray}
Then, go to Step 1) with $i$ replaced by $i+1$.
\end{enumerate}


\tikzstyle{decision} = [diamond, draw, 
    text width=6.3em, text badly centered, node distance=2.3cm, inner sep=0pt, aspect=3.5]
\tikzstyle{block} = [rectangle, draw, 
    text width=3.5cm, rounded corners, minimum height=3em]
\tikzstyle{wideblock} = [rectangle, draw, 
    text width=7.5cm, text centered, rounded corners, minimum height=3em]
        \tikzstyle{wideblock2} = [rectangle, draw, 
    text width=6.5cm, text centered, rounded corners, minimum height=3em]
\tikzstyle{line} = [draw, -latex']
\tikzstyle{cloud} = [draw, ellipse,fill=blue!10, node distance=3cm,
    minimum height=2em]
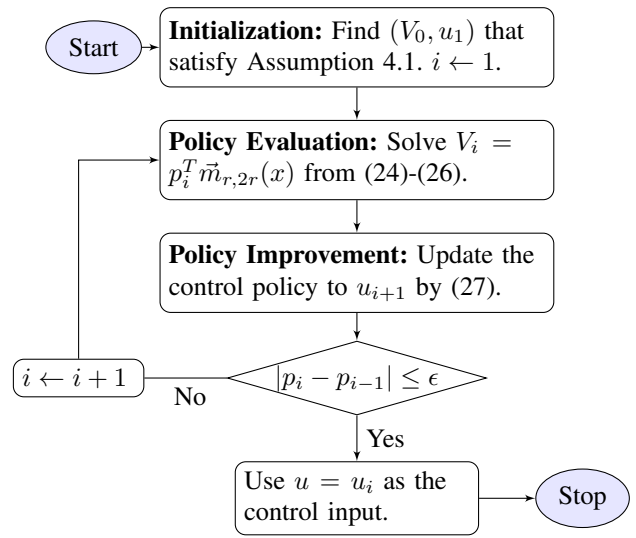
\begin{figure}[!b]
\begin{tikzpicture}[node distance = 3.5cm, auto]
    \node [cloud, node distance = 2cm] (start) {Start};
    \node [block, right of=start, text width=5cm] (init) {{\bf Initialization:}  Find $(V_0, u_1)$ that satisfy Assumption \ref{ass:polyall}. $i\gets 1$.};
    %
    %
    \node [block, below of=init, text width=5cm, node distance=1.5cm] (clct) {{\bf Policy Evaluation:} Solve $V_i = p_i^T\vec{m}_{r,2r}(x)$ from \eqref{eq:sosopt1}-\eqref{eq:sosopt4}.};
    %
    \node[block, below of=clct, text width=5cm, node distance=1.5cm](main)
    {{\bf Policy Improvement:} Update the control policy to $u_{i+1}$ by \eqref{eq:sosu}.};
    \node [decision, below of=main, node distance=1.4cm] (decide) {$|p_i-p_{i-1}|\le \epsilon$};
    \node [block, below of=decide, node distance=1.6cm, text width=3cm, minimum height=.5cm] (realcontrol) {Use $u = u_{i}$ as the control input.};
    \node [cloud, right of=realcontrol, node distance=3cm] (stop) {Stop};
    \node [block, left of=decide, minimum height=.5cm, text width=1.5cm, node distance=3.7cm](kp1){$i\gets i+1$};
    \path [line] (start) -- (init);
    \path [line] (init) -- (clct);
    \path [line] (clct) -- (main);
    \path [line] (main) -- (decide);
    \path [line] (decide) -- node {Yes} (realcontrol);
    \path [line] (decide)  --node {~No}(kp1)  |- (clct);
    \path [line] (realcontrol) --(stop);
  \end{tikzpicture}
  \caption{Flowchart of the SOS-programming-based policy iteration.}
  \label{flowchart1}
\end{figure}

%
%


A flowchart of the practical implementation is given in Figure \ref{flowchart1}, where $\epsilon>0$ is a predefined threshold to balance the exploration/exploitation trade-off. In practice, a larger $\epsilon$ may lead to shorter exploration time and therefore will allow the system to implement the control policy and terminate the exploration noise sooner. On the other hand, using a smaller $\epsilon$ allows the learning system to better improve the control policy but needs longer learning time.

\begin{remark}
The optimization problem \eqref{eq:sosopt1}-\eqref{eq:sosopt4} is a well-defined SOS program \cite{blekherman2013semidefinite}. Indeed, the  objective function \eqref{eq:sosopt1} is linear in $p$, since for any $V=p^T\vec{m}_{2,2r}(x)$, we have $\int_{\Omega}V(x)dx=c^Tp$, with $c=\int_{\Omega}\vec{m}_{2,2r}(x)dx$.
\end{remark}

\begin{thm}
 \label{cor:cor4sos}
 Under Assumptions \ref{ass:HJB} and  \ref{ass:polyall}, the following are true, for $i=1,2,\cdots$.
\begin{enumerate}
\item[1)]
The SOS program \eqref{eq:sosopt1}-\eqref{eq:sosopt4} has a nonempty feasible set.

\item[2)] The closed-loop system comprised of \eqref{eq:sys1} and $u=u_{i}$
is globally asymptotically stable at the origin.
\item[3)] $V_{i}\in \mathcal{P}$. In particular, the following inequalities hold:
\begin{eqnarray}
V^{\rm o}(x_0)\le V_{i}(x_{0})\le V_{i-1}(x_{0}),~~\forall x_{0}\in\mathbb{R}^{n}.
\end{eqnarray}
\item[4)] There exists $V^{*}(x)$ satisfying $V^{*}(x)\in\mathbb{R}[x]_{2,2r}\cap \mathcal{P}$, such that, for any $x_0\in\mathbb{R}^n$,  $\lim\limits _{i\rightarrow\infty}V_{i}(x_0) =  V^{*}(x_0)$.

\item[5)] Along the solutions of the system comprised of \eqref{eq:sys1} and \eqref{eq:optimality},
the following inequalities hold:
\begin{eqnarray}
0\le V^*(x_0)-V^{\rm o}(x_0)\le -\int_0^{\infty}\mathcal{H}(V^*(x(t)))dt.
\end{eqnarray}
\end{enumerate}
\end{thm}

\begin{IEEEproof}
1) We prove by mathematical induction.

i) Suppose $i=1$, under Assumption \ref{ass:polyall}, we know $\mathcal{L}(V_0,u_1)$ is SOS. Hence, $V=V_0$ is a feasible solution to the problem \eqref{eq:sosopt1}-\eqref{eq:sosopt4}.

ii) Let $V=V_{j-1}$ be an optimal solution to the problem \eqref{eq:sosopt1}-\eqref{eq:sosopt4} with $i=j-1>1$. We show that $V=V_{j-1}$ is a feasible solution to the same problem with $i=j$.

Indeed, by definition,
\[u_{j}=-\frac{1}{2}R^{-1}g^T\nabla V_{j-1},\]
  and
\begin{eqnarray*}
\mathcal{L}(V_{j-1},u_j)&=&-\nabla V_{j-1}^T(f+gu_j)-r(x,u_j)\\
&=&\mathcal{L}(V_{j-1},u_{j-1})-\nabla V_{j-1}^T g(u_j-u_{j-1})\\
&& +u_{j-1}^TRu_{j-1}-u_{j}^TRu_{j}\\
&=&\mathcal{L}(V_{j-1},u_{j-1})+|u_j-u_{j-1}|_R^2.
\end{eqnarray*}

Under the induction assumption, we know $V_{j-1}\in\mathbb{R}[x]_{2,2r}$ and $\mathcal{L}(V_{j-1},u_{j-1})$ is SOS. Hence, $\mathcal{L}(V_{j-1},u_{j})$ is SOS. As a result, $V_{j-1}$ is a feasible solution to the SOS program \eqref{eq:sosopt1}-\eqref{eq:sosopt4} with $i=j$.

2) Again, we prove by induction.

i) Suppose $i=1$, under Assumption \ref{ass:polyall}, $u_1$ is globally stabilizing. Also, we can show that $V_1\in\mathcal{P}$. Indeed, for each $x_0\in\mathcal{R}^n$ with $x_0\neq 0$, we have
\begin{eqnarray}
V_1(x_0)\ge\int_0^{\infty} r(x,u_1)dt>0. \label{eq:v1ger}
\end{eqnarray}
By \eqref{eq:v1ger} and the constraint \eqref{eq:sosopt4}, under Assumption \ref{ass:HJB} it follows that
\begin{eqnarray}
V^{\rm o}\le V_1\le V_0.  \label{eq:proper}
\end{eqnarray}
Since both $V^{\rm o}$ and $V_0$ are assumed to belong to $\mathcal{P}$, we conclude that $V_1\in\mathcal{P}$.

ii) Suppose $u_{i-1}$ is globally stabilizing, and $V_{i-1}\in\mathcal{P}$ for $i>1$. Let us show that $u_{i}$ is globally stabilizing, and $V_{i}\in\mathcal{P}$.

Indeed, along the solutions of the closed-loop system composed of \eqref{eq:sys1} and $u=u_i$, it follows that
\begin{eqnarray*}
\dot{V}_{i-1} & = & \nabla V_{i-1}^T(f+gu_{i}) \\
              & = &-\mathcal{L}(V_{i-1},u_i)-r(x,u_i)\\
               &\le&-q(x).
\end{eqnarray*}

Therefore, $u_i$ is globally stabilizing, since $V_{i-1}$ is a well-defined Lyapunov function for the system. In addition, we have
\begin{eqnarray}
V_i(x_0)\ge\int_0^{\infty} r(x,u_i)dt>0,~~\forall x_0\neq 0.
\end{eqnarray}
Similarly as in \eqref{eq:proper}, we can show
\begin{eqnarray}
V^{\rm o}(x_0)\le V_i(x_0)\le V_{i-1}(x_0),~~\forall x_0\in\mathbb{R}^n, \label{eq:general_inequality}
\end{eqnarray}
and conclude that $V_i\in\mathcal{P}$.

3) The two inequalities have been proved in \eqref{eq:general_inequality}.

4) By 3), for each $x\in\mathbb{R}^{n}$, the sequence $\{V_{i}(x)\}_{i=0}^{\infty}$
is monotonically decreasing with $0$ as its lower bound. Therefore, the
limit exists, i.e., there exists $V^{*}(x)$,
such that $\lim\limits _{i\rightarrow\infty}V_{i}(x)=V^{*}(x)$. Let $\{p_i\}_{i=1}^\infty$ be the sequence such that $V_i=p_i^T\vec{m}_{2,2r}(x)$. Then, we know $\lim\limits_{i\rightarrow\infty}p_i=p^*\in\mathbb{R}^{n_{2r}}$, and therefore $V^*=p^{*T}\vec{m}_{2,2r}(x)$. Also, it is easy to show $V^{\rm o}\le V^*\le V_0$. Hence, $V^*\in\mathbb{R}[x]_{2,2r}\cap \mathcal{P}$.
%

5) Let
\begin{eqnarray}
u^*=-\frac{1}{2}R^{-1}g^T\nabla V^*,\label{eq:ustar}
\end{eqnarray}
Then, by 4) we know
\begin{eqnarray}
\mathcal{H}(V^*)=-\mathcal{L}(V^*,u^*)\le 0,
\end{eqnarray}
which implies that $V^*$ is a feasible solution to Problem \ref{prob:relaxed_HJB}.
The inequalities in 5) can thus be obtained by the fourth property in Theorem \ref{thm:relaxed_HJB}.

The proof is thus complete.
%
\end{IEEEproof}

\begin{remark}
Notice that Assumption \ref{ass:polyall}  holds for any controllable linear systems. For general nonlinear systems, such a pair $(V_0, u_1)$ satisfying Assumption 4.1 may not always exist, and in this paper we focus on polynomial systems that satisfy Assumption \ref{ass:polyall}. For this class of systems,  we assume that both $V_0$ and $u_1$  have to be determined before executing the proposed algorithm. The search of $(V_0, u_1)$ is not trivial, because it amounts to solving some bilinear matrix inequalities (BMI) \cite{prajna2004nonlinear}. However, this problem has been actively studied in recent years, and several applicable approaches have been developed. For example, a Lyapunov based approach utilizing state-dependent linear matrix inequalities has been studied in \cite{prajna2004nonlinear}. This method has been generalized to uncertain nonlinear polynomial systems in \cite{xu2009simultaneous} and \cite{huang2013robust}. In \cite{Horowitz_ACC2014} and \cite{Horowitz_Preprint2014}, the authors proposed a solution for a stochastic HJB. It is shown that this method gives a control Lyapunov function for the deterministic system.
\end{remark}

\begin{remark}
\label{rm:nonPolyU}
Notice that the control policies considered in the proposed algorithm can be extended to polynomial fractions. Indeed, instead of requiring  $\mathcal{L}(V_0,u_1)$ to be SOS in Assumption \ref{ass:polyall}, let us assume
\begin{eqnarray}
\alpha^2(x)\mathcal{L}(V_0,u_{1}) {\rm~is~SOS}  \label{eq:sos_newcs0}
\end{eqnarray}
where $\alpha(x)>0$ is an arbitrary  polynomial. Then, the initial control policy $u_1$ can take the form of $\displaystyle u _1 = {\alpha(x)}^{-1}v(x)$ , with $v(x)$ a column vector of polynomials. Then, we can relax the constraint \eqref{eq:sosopt3} to the following:
\begin{eqnarray}
\alpha^2(x)\mathcal{L}(V,u_{i}) {\rm~is~SOS}  \label{eq:sos_newcs}
\end{eqnarray}
and it is easy to see the SOS-based policy iteration algorithm can still proceed with this new constraint.
\end{remark}

\begin{remark}
\label{rm:nonConstantR}
In addition to Remark \ref{rm:nonPolyU}, if $R(x)$ is not constant, by definition, there exists a symmetric matrix of polynomials $R^*(x)$, such that
\begin{eqnarray}
R(x)R^*(x) = {\rm det}(R(x)) I_{m}
\end{eqnarray}
with ${\rm det}(R(x))>0$.

As a result, the policy improvement step \eqref{eq:sosu} gives
\begin{eqnarray}
u_{i+1} = -\frac{1}{2{\rm det}(R(x)) }R^*(x)g^T(x)\nabla V_{i}(x).
\end{eqnarray}
which is a vector of polynomial fractions.

Now, if we select $\alpha(x)$ in \eqref{eq:sos_newcs} such that ${\rm det}(R(x))$ divides  $\alpha(x)$, we see $\alpha^2(x)\mathcal{L}(V,u_{i})$ is a polynomial, and the proposed SOS-based policy iteration algorithm can still proceed with \eqref{eq:sosopt3} replaced by \eqref{eq:sos_newcs}. Following the same reasoning as in the proof of Theorem \ref{cor:cor4sos}, it is easy to show that the solvability and the convergence properties of the proposed policy iteration algorithm will not be affected with the relaxed constraint \eqref{eq:sos_newcs}.
\end{remark}

\section{Global adaptive dynamic programming for uncertain polynomial systems}
\label{sec:online_poly}
The proposed policy iteration method requires the perfect knowledge of $f$ and $g$. In practice, precise system knowledge may be difficult to obtain. Hence, in this section, we develop an online learning method based on the idea of ADP to implement the iterative scheme with real-time data, instead of identifying first the system dynamics. For simplicity, here again we assume $R(x)$ is a constant matrix. The results can be straightforwardly extended to non-constant $R(x)$,  using the same idea as discussed in Remark \ref{rm:nonConstantR}.


To begin with, consider the system
\begin{eqnarray}
\dot x = f+g(u_i+e) \label{eq:sys_with_exploration_noise}
\end{eqnarray}
where $u_i$ is a feedback control policy and $e$ is a bounded time-varying function, known as the exploration noise, added for the learning purpose.

\subsection{Forward completeness}

Due to the existence of the exploration noise, we need to first show if solutions of system \eqref{eq:sys_with_exploration_noise} exist globally, for positive time. For this purpose, we will show that the system \eqref{eq:sys_with_exploration_noise} is forward complete \cite{angeli1999forward}.

\begin{definition}[\cite{angeli1999forward}]
Consider system \eqref{eq:sys_with_exploration_noise} with $e$ as the input. The system is called {\it forward complete} if, for any initial condition $x_0\in\mathbb{R}^n$ and every input signal $e$, the corresponding solution of system \eqref{eq:sys_with_exploration_noise} is defined for all $t\ge0$.
\end{definition}

\begin{lemma}
\label{le:FC}
Consider system \eqref{eq:sys_with_exploration_noise}.
Suppose $u_i$ is a globally stabilizing control policy and there exists $V_{i-1}\in \mathcal{P}$, such that $\nabla V_{i-1}(f+gu_i)+u_i^TRu_i\le0$. Then, the system \eqref{eq:sys_with_exploration_noise} is forward complete.
\end{lemma}

\begin{IEEEproof} Under Assumptions \ref{ass:HJB} and \ref{ass:polyall}, by  Theorem \ref{cor:cor4sos} we know $V_{i-1}\in\mathcal{P}$. Then, by completing the squares, it follows that
\begin{eqnarray*}
\nabla V_{i-1}^T(f+gu_i+ge)
&\le& -u_i^TRu_i-2u_{i}^TRe\\
&=& -|u_i+e|^2_R+|e|^2_R\\
&\le& |e|^2_R\\
&\le& |e|^2_R+V_{i-1}.
\end{eqnarray*}
According to \cite[Corollary 2.11]{angeli1999forward}, the system \eqref{eq:sys_with_exploration_noise} is forward complete.
\end{IEEEproof}

By Lemma \ref{le:FC} and Theorem \ref{cor:cor4sos}, we immediately have the following Proposition.
\begin{proposition}
Under Assumptions \ref{ass:HJB} and \ref{ass:polyall}, let $u_i$ be a feedback control policy obtained at the $i$-th iteration step in the proposed policy iteration algorithm \eqref{eq:sosopt1}-\eqref{eq:sosu} and $e$ be a bounded time-varying function. Then, the closed-loop system \eqref{eq:sys1} with $u=u_i+e$ is forward complete.
\end{proposition}

\subsection{Online implementation}

Under Assumption \ref{ass:polyall}, in the SOS-based policy iteration,  we have
\begin{eqnarray}
\mathcal{L}(V_i,u_i)\in\mathbb{R}[x]_{2,2d},~~\forall i>1,
\end{eqnarray}
if the integer $d$ satisfies
\begin{eqnarray*}
d &\ge& \frac{1}{2}{\rm max} \left\{{\rm deg}(f) + 2r-1, {\rm deg}(g) +2(2r-1) ,  \right.\\
&&\left. {\rm deg}(Q), 2(2r-1) + 2{\rm deg}(g)\right\}.
\end{eqnarray*}
Also,  $u_i$ obtained from the proposed policy iteration algorithm satisfies $u_i\in\mathbb{R}[x]_{1,d}$.

Hence, there exists a constant matrix $K_{i}\in\mathbb{R}^{m\times n_{d}}$, with $n_d=(^{n+d}_{~d})-1$, such that $u_i=K_i\vec{m}_{1,d}(x)$. Also, suppose there exists a constant vector $p\in\mathbb{R}^{ n_{2r}}$, with $n_{2r}=(^{n+2r}_{~2r})-n-1$,  such that $V=p^T\vec{m}_{2,2r}(x)$. Then, along the solutions  of the system \eqref{eq:sys_with_exploration_noise}, it follows that
\begin{eqnarray}
\dot V&=&\nabla V^T\left(f+gu_i\right)+\nabla V^Tge\nonumber\\
&=&-r(x,u_i)-\mathcal{L}(V,u_i)+\nabla V^Tge\nonumber\\
&=&-r(x,u_i)-\mathcal{L}(V,u_i) \nonumber\\
&&+ (R^{-1}g^T\nabla V)^T Re  \label{eq:online_polyn0}
\end{eqnarray}
Notice that the two terms $\mathcal{L}(V,u_i)$ and $R^{-1}g^T\nabla V$ rely on $f$ and $g$. Here, we are interested in solving them without identifying $f$ and $g$.

To this end, notice that for the same pair $(V,u_i)$ defined above,  we can find a constant vector $l_p\in\mathbb{R}^{n_{2d}}$, with $n_{2d}=(^{n+2d}_{~2d})-d-1$, and a constant matrix $K_p\in\mathbb{R}^{m\times n_{d}}$, such that
\begin{eqnarray}
\mathcal{L}(V, u_i) &=& l_p^T\vec{m}_{2,2d}(x) \label{eq:iota}\\
-\frac{1}{2}R^{-1}g^T\nabla V &=&K_p\vec{m}_{1,d}(x) \label{eq:kappa}
\end{eqnarray}
Therefore, calculating $\mathcal{L}(V,u_i)$ and $R^{-1}g^T\nabla V$  amounts to finding $l_p$ and $K_p$.

%


Substituting \eqref{eq:iota} and \eqref{eq:kappa} in \eqref{eq:online_polyn0}, we have
\begin{eqnarray}
\dot V =-r(x,u_i)-l_p^T\vec{m}_{2,2d}(x)-2\vec{m}_{1,d}^T(x)K_p^TRe \label{eq:online_polyn}
\end{eqnarray}

Now, integrating the terms in \eqref{eq:online_polyn} over the interval $[t,t+\delta t]$, we have
\begin{eqnarray}
&&p^T\left[\vec{m}_{2,2r}(x(t))-\vec{m}_{2,2r}(x(t+\delta t))\right]\nonumber\\
&=& \int_{t}^{t+\delta t}\left(r(x,u_{i})+l_p^T\vec{m}_{2,2d}(x)\right.\nonumber\\
&& \left.+2\vec{m}_{1,d}^T(x)K_p^TRe\right)dt \label{eq:intADPpoly}
\end{eqnarray}

Eq. \eqref{eq:intADPpoly} implies that, $l_p$ and $K_p$ can be directly calculated by using real-time online data, without knowing the precise knowledge of $f$ and $g$.

To see how it works, let us define the following matrices: $\sigma_e\in\mathbb{R}^{n_{2d}+mn_{d}}$, $\Phi_{i}\in\mathbb{R}^{q_i\times (n_{2d}+mn_d)}$, $\Xi_i\in\mathbb{R}^{q_{i}}$, $\Theta_i\in\mathbb{R}^{q_{i}\times n_{2r}}$, such that
\setlength{\arraycolsep}{5pt}
\begin{eqnarray*}
\sigma_e&=&-\left[\begin{array}{cc}\vec{m}_{2,2d}^T&2\vec{m}_{1,d}^T\otimes e^{T}R\end{array}\right]^T,\\
\Phi_{i}&=&
\left[\begin{array}{cccc}
\int_{t_{0,i}}^{t_{1,i}}\sigma_e dt&
\int_{t_{1,i}}^{t_{2,i}}\sigma_edt&
\cdots&
\int_{t_{q_i-1,i}}^{t_{q_i,i}}\sigma_edt
\end{array}\right]^T,\\
\Xi_{i}&=&\left[\begin{array}{ccc}
\int_{t_{0,i}}^{t_{1,i}}r(x,u_{i})dt&
\int_{t_{1,i}}^{t_{2,i}}r(x,u_{i})dt&
\cdots\end{array}\right.\\
&&\left.\begin{array}{c}\int_{t_{q_{i}-1,i}}^{t_{q_{i},i}}r(x,u_{i})dt\end{array}
\right]^T,\\
\Theta_{i}&=&\left[\begin{array}{cccc}
\vec{m}_{2,2r}|_{t_{0,i}}^{t_{1,i}}&
\vec{m}_{2,2r}|_{t_{1,i}}^{t_{2,i}}&
\cdots&
\vec{m}_{2,2r}|_{t_{q_{i}-1,i}}^{t_{q_{i},i}}
\end{array}\right]^T.
\end{eqnarray*}
\setlength{\arraycolsep}{5pt}

Then, \eqref{eq:intADPpoly} implies
\begin{eqnarray}
\Phi_{i}\left[\begin{array}{c} l_p\\{\rm vec}(K_p)\end{array}\right]=\Xi_i+\Theta_i p. \label{eq:online0}
\end{eqnarray}

Notice that any pair of $(l_p, K_p)$ satisfying \eqref{eq:online0} will satisfy the constraint between $l_p$ and $K_p$ as implicitly indicated in \eqref{eq:online_polyn} and \eqref{eq:intADPpoly}.

\begin{assume}
\label{ass:rank_condition}
For each $i=1,2,\cdots$, there exists an integer $q_{i0}$, such that
 the following rank condition holds,
\begin{eqnarray}
{\rm rank}(\Phi_{i})=n_{2d}+mn_{d}, \label{eq:rank1}
\end{eqnarray}
if $q_i\ge q_{i0}$.
\end{assume}

\begin{remark}
This rank condition \eqref{eq:rank1} is in the spirit of persistency of excitation (PE) in adaptive control ({\it e.g.} \cite{Ioannou1996adaptive, tao2003adaptive}) and is a necessary condition for parameter convergence.
\end{remark}

Given $p\in\mathbb{R}^{n_{2r}}$ and $K_i\in\mathbb{R}^{m\times {n_{d}}}$, suppose Assumption \ref{ass:rank_condition} is satisfied and $q_i\ge q_{i0}$ for all $i=1,2, \cdots$. Then, it is easy to see that the values of  $l_p$ and $K_p$ can be uniquely determined from \eqref{eq:intADPpoly}. Indeed,
\begin{eqnarray}\left[\begin{array}{c}l_p\\{\rm
vec}(K_p)\end{array}\right]=\left(\Phi_{i}^T\Phi_{i}\right)^{-1}\Phi_{i}^T\left(\Xi_i+\Theta_ip\right)
\end{eqnarray}


Now, the ADP-based online learning method is given below, and a flowchart is provided in Figure \ref{flowchart2}.
\begin{enumerate}
\item[1)] {\it Initialization:} \\
Fine the pair $(V_0, u_1)$ that satisfy Assumption \ref{ass:polyall}. Let $p_0$ be the constant vector such that $V_0=p_0^T\vec{m}_{2,2r}(x)$, Let $i=1$.
\item[2)] {\it Collect online data:}\\
Apply $u=u_i+e$ to the system and compute the data matrices $\Phi_i$, $\Xi_i$, and $\Theta_i$, until $\Phi_i$ is of full column rank.
\item[3)] \textit{Policy evaluation and improvement:}\\
Find an optimal solution $(p_i, K_{i+1})$ to the following SOS program
\begin{eqnarray}
\min_{p,K_p}~c^Tp &&\label{eq:onlinesdpopt1}\\
{\rm ~s.t.}~~
\Phi_{i}\left[\begin{array}{c} l_p \\{\rm vec}(K_p)\end{array}\right]&=&\Xi_i+\Theta_i p \label{eq:onlinesdpopt2}\\
l_p^T \vec{m}_{2,2d}(x) &{\rm is}& {\rm SOS}\\
 (p_{i-1}-p)^T\vec{m}_{2,2r}(x) &{\rm is}& {\rm SOS} \label{eq:onlinesdpopt4}
\end{eqnarray}
where $c = \int_{\Omega} \vec{m}_{2,2r}(x)dx$.

Then, denote $V_i=p_i^T\vec{m}_{2,2r}(x)$, $u_{i+1}=K_{i+1}\vec{m}_{1,d}(x)$, and go to Step 2) with $i\gets i+1$.
\end{enumerate}

\begin{thm} Under Assumptions \ref{ass:controlability},  \ref{ass:polyall} and \ref{ass:rank_condition}, the following properties hold.
\begin{enumerate}
\item The optimization problem \eqref{eq:onlinesdpopt1}-\eqref{eq:onlinesdpopt4} has a nonempty feasible set.
\item The sequences $\{V_i\}_{i=1}^\infty$ and $\{u_i\}_{i=1}^\infty$ satisfy the properties 2)-5) in Theorem \ref{cor:cor4sos}.
\end{enumerate}
\end{thm}

\begin{IEEEproof}
Given $p_i\in\mathbb{R}^{n_{2r}}$, there exists a constant matrix $K_{p_i}\in\mathbb{R}^{m\times n_{d}}$ such that $(p_i, K_{p_i})$ is a feasible solution to the optimization problem \eqref{eq:onlinesdpopt1}-\eqref{eq:onlinesdpopt4} if and only if $p_i$ is a feasible solution to the SOS program \eqref{eq:sosopt1}-\eqref{eq:sosopt4}. Therefore, by Theorem \ref{cor:cor4sos}, 1) holds. In addition, since the two optimization problems share the identical objective function, we know that if $(p_i, K_{p_i})$ is a feasible solution to the optimization problem \eqref{eq:onlinesdpopt1}-\eqref{eq:onlinesdpopt4}, $p_i$ is also an optimal solution to the SOS program \eqref{eq:sosopt1}-\eqref{eq:sosopt4}. Hence, the theorem can be obtained from Theorem \ref{cor:cor4sos}.
\end{IEEEproof}

%
%
\tikzstyle{decision} = [diamond, draw, 
    text width=6.3em, text badly centered, node distance=2.3cm, inner sep=0pt, aspect=3.5]
\tikzstyle{block} = [rectangle, draw, 
    text width=3.5cm, rounded corners, minimum height=3em]
\tikzstyle{wideblock} = [rectangle, draw, 
    text width=7.5cm, text centered, rounded corners, minimum height=3em]
        \tikzstyle{wideblock2} = [rectangle, draw, 
    text width=6.5cm, text centered, rounded corners, minimum height=3em]
\tikzstyle{line} = [draw, -latex']
\tikzstyle{cloud} = [draw, ellipse,fill=blue!10, node distance=3cm,
    minimum height=2em]
\begin{figure}[!b]
\begin{tikzpicture}[node distance = 1.5cm, auto]
    \node [cloud, node distance = 1.5cm] (start) {Start};
    \node [block, below of=start, text width=6cm, node distance = 1.2cm] (init) {{\bf Initialization:}  Find $(V_0, u_1)$ that satisfy Assumption \ref{ass:polyall}. $i\gets 1$.};
    \node [block, below of=init, text width=6cm, node distance = 1.8cm] (oldata) {{\bf Collect Online Data:}  Apply $u_i$ as the control input. Compute the matrices $\Phi_i$, $\Xi_i$, and $\Theta_i$ until Assumption \ref{ass:rank_condition} is satisfied.};
    %
    %
    %
    \node [block, below of=oldata, text width=6cm, node distance=2.2cm] (clct) {{\bf Policy Evaluation and Improvement:}  Obtain $p_i$ and $K_{i+1}$ by solving \eqref{eq:onlinesdpopt1}-\eqref{eq:onlinesdpopt4}. Let $V_i = p_i^T\vec{m}_{2,2r}(x)$ and $u_{i+1}=K_{i+1}\vec{m}_{1,d}(x).$};
    %
    \node [decision, below of=clct, node distance=1.7cm] (decide) {$|p_i-p_{i-1}|\le \epsilon$};
    \node [block, below of=decide, node distance=1.6cm, text width=3cm, minimum height=.5cm] (realcontrol) {Apply $u = u_{i}$ as the control input.};
    \node [cloud, right of=realcontrol, node distance=3cm] (stop) {Stop};
    \node [block, left of=decide, minimum height=.5cm, text width=1.5cm, node distance=3.7cm](kp1){$i\gets i+1$};
    \path [line] (start) -- (init);
    \path [line] (init) -- (oldata);
    \path [line] (oldata) -- (clct);
    \path [line] (clct) -- (decide);
    \path [line] (decide) -- node {Yes} (realcontrol);
    \path [line] (decide)  --node {~No}(kp1)  |- (oldata);
    \path [line] (realcontrol) --(stop);
  \end{tikzpicture}
  \caption{Flowchart of ADP-based online control method.}
  \label{flowchart2}
\end{figure}
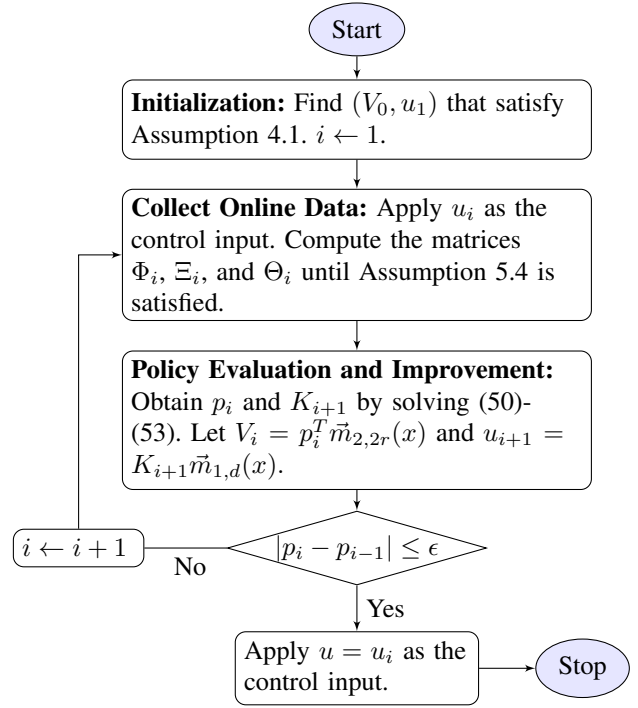

\begin{remark}
By Assumption \ref{ass:polyall}, $u_1$ must be a globally stabilizing control policy. Indeed, if it is only locally stabilizing, there are two reasons why the algorithm may not proceed. First, with a locally stabilizing control law, there is no way to guarantee the {\it forward completeness} of solutions, or avoid finite escape, during the learning phase. Second,  the region of attraction associated with the new control policy is not guaranteed to be global. Such a counterexample is $\dot x = \theta x^3 +u$ with unknown positive $\theta$, for which the choice of a locally stabilizing controller $u_1= -x$ does not lead to a globally stabilizing suboptimal controller.
\end{remark}

\section{Applications}


\subsection{A scalar nonlinear polynomial system}

Consider the following polynomial system
\begin{eqnarray}
\dot x = a x^2+bu \label{eq:scalarsys}
\end{eqnarray}
where $x\in\mathbb{R}$ is the system state, $u\in\mathbb{R}$ is the control input, $a$ and $b$, satisfying $a\in [0,0.05]$ and $b\in [0.5,1]$, are uncertain constants. The cost to be minimized is defined as
$J(x_0,u) = \int_0^{\infty} (0.01 x^2+ 0.01 x^4 +u^2)dt$.
%
An initial stabilizing control policy can be selected as $u_1=-0.1x-0.1x^3$, which globally asymptotically stabilizes system \eqref{eq:scalarsys},
for any $a$ and $b$ satisfying the given range.    Further, it is easy to see that $V_0=10(x^2+x^4)$ and $u_1$ satisfy Assumption \ref{ass:polyall} with $r=2$. In addition, we set $d=3$.

Only for the purpose of simulation, we set $a=0.01$, $b=1$, and $x(0)=2$. $\Omega$ is specified as $\Omega=\{x|x\in\mathbb{R}~{\rm and}~|x|\le 1\}$. The proposed global ADP method is applied with the control policy updated after every five seconds, and convergence is attained after five iterations, when $|p_{i}-p_{i-1}|\le 10^{-3}$.
Hence, the constant vector $c$ in the objective function \eqref{eq:onlinesdpopt1} is computed as $c=[\begin{array}{ccc}\frac{2}{3}& 0& \frac{2}{5}\end{array}]^T$. The exploration noise is set to be $e=0.01(\sin(10t)+\sin(3t)+\sin(100t))$, which is turned off after the fourth iteration.

\begin{figure}
\centering
\includegraphics[width=.4\textwidth]{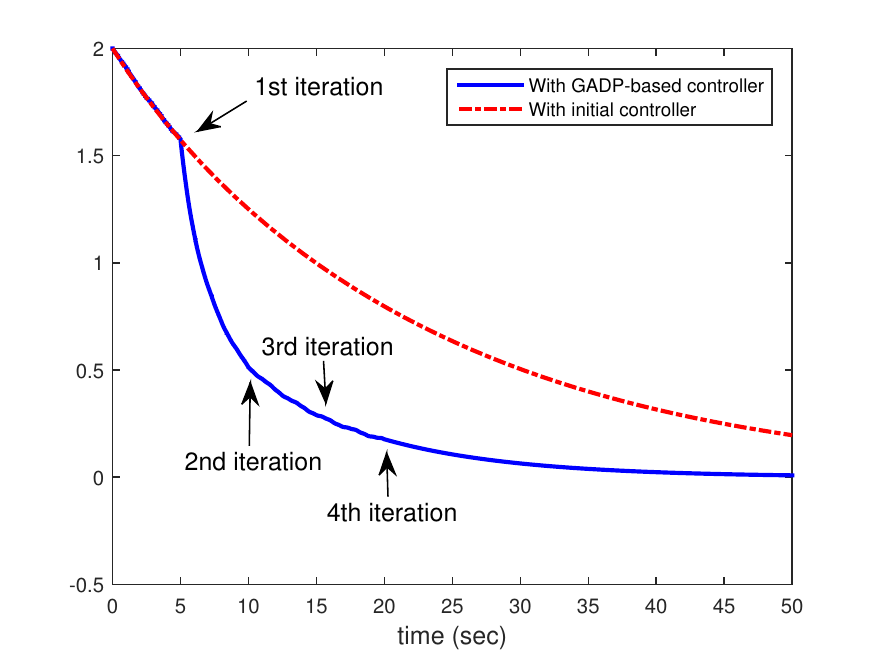}
\caption{System state trajectory.}
\label{fig_ex1_state}
\end{figure}

\begin{figure}
\centering
\includegraphics[width=0.4\textwidth]{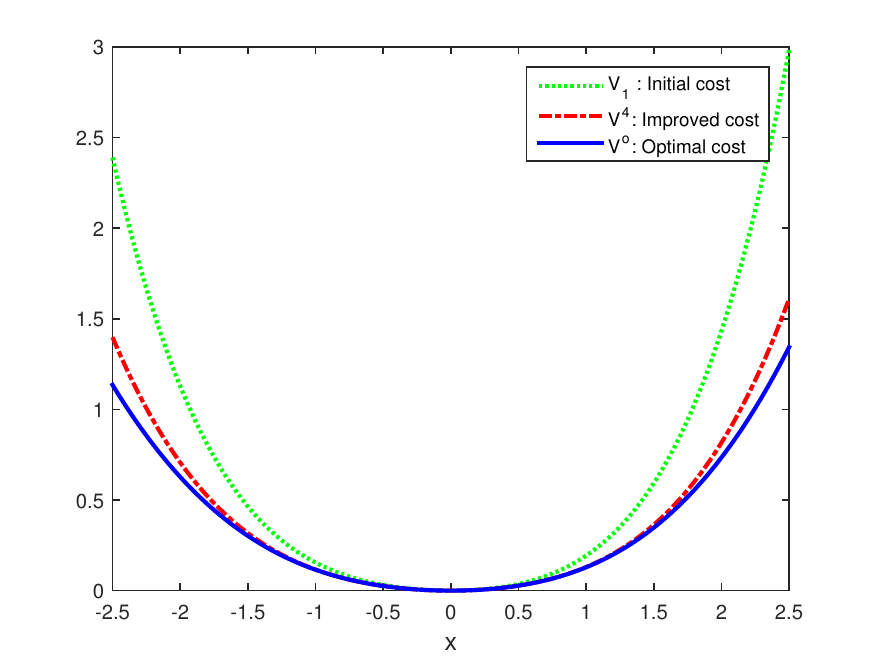}
\caption{Comparison of the value functions.}
\label{fig_ex1_cost}
\end{figure}

\begin{figure}[!t]
\centering
\includegraphics[width=0.4\textwidth]{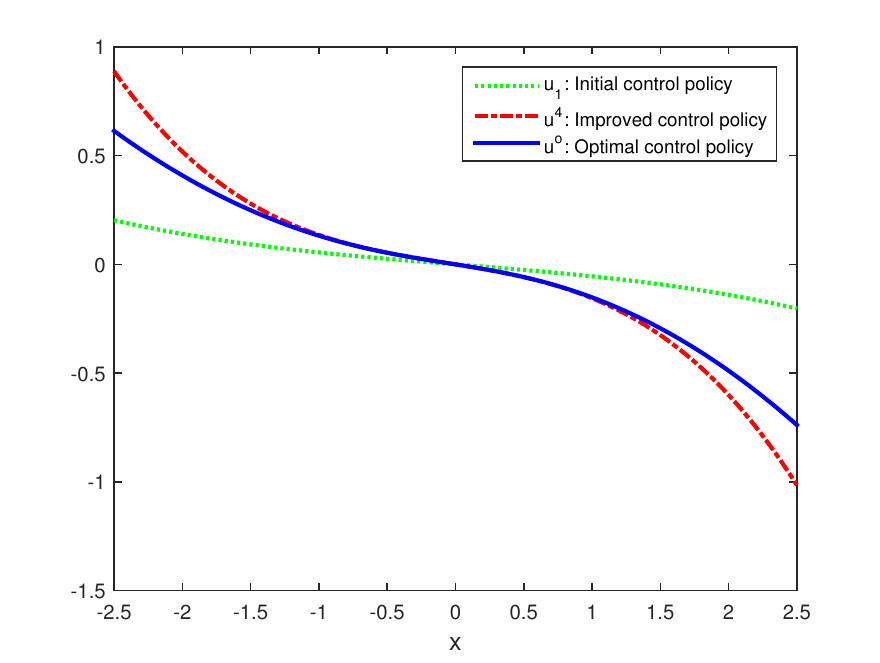}
\caption{Comparison of the control policies.}
\label{fig_ex1_control_compare}
\end{figure}

The simulated state trajectory is shown in Figure \ref{fig_ex1_state}, where the control policy is updated every five seconds until convergence is attained. The suboptimal control policy and the cost function obtained after four iterations are $V^*=0.1020x^2+0.007x^3+0.0210x^4$ and $u^*=-0.2039x-0.02x^2-0.0829x^3$.
For comparison purpose, the exact optimal cost and the control policy are given below.
$V^{\rm o} = \frac{x^3}{150} + \frac{(\sqrt{101x^2 + 100})^3}{15150} - \frac{20}{303}$ and $u^{\rm o} = -\frac{x^2\sqrt{101x^2 + 100} + 101x^4 + 100x^2}{100\sqrt{101x^2 + 100}}$. Figures \ref{fig_ex1_cost} and \ref{fig_ex1_control_compare} show the comparison of the suboptimal control policy with respect to the exact optimal control policy and the initial control policy.



\subsection{Fault-tolerant control}

Consider the following example  \cite{xu2009simultaneous}
\begin{eqnarray}
\left[\begin{array}{c}\dot x_1\\ \dot x_2\end{array}\right]
\!\!=\!\!
\left[\begin{array}{c} -x_1^3-x_1x_2^2+x_1x_2 \\x_1 +2x_2\end{array}\right]
\!\!+\!\! 
\left[\begin{array}{cc}
0&\beta_1\\
\beta_2&\beta_1
\end{array}\right]u 
\label{eq:sysjet}
\end{eqnarray}
where $\beta_1,\beta_2\in[0.5,1]$ are the uncertain parameters. This system can be considered as having a loss-of-effectiveness fault \cite{xu2009simultaneous}, since the actuator gains are smaller than the commanded position ($\beta_1=\beta_2=1$).

Using SOS-related techniques, it has been shown in \cite{xu2009simultaneous} that the following robust control policy  can globally asymptotically stabilize the system \eqref{eq:sysjet} at the origin.
\begin{eqnarray}
u_1 = \left[\begin{array}{c}u_{1,1} \\
u_{1,2}
\end{array}
\right]=\left[
\begin{array}{c}
 10.283x_1 -13.769x_2\\
 -10.7x_1 -3.805x_2
\end{array}
\right] \label{eq:ex2_u1}
\end{eqnarray}
However, optimality of the closed-loop system has not been well addressed.

Our control objective is to improve the control policy using the proposed ADP method such that we can reduce the following cost
\begin{eqnarray}
J(x_0,u)=\int_0^{\infty}(x_1^2+x_2^2+u^Tu) dt
\end{eqnarray}
and at the same time guarantee global asymptotic stability. In particular, we are interested to improve the performance of the closed-loop system in the set $\Omega=\{x|x\in\mathbb{R}^2~{\rm and}~|x|\le 1\}$.

By solving the following feasibility problem using SOSTOOLS \cite{sostools, prajna2002introducing}
\begin{eqnarray}
&V\in\mathbb{R}[x]_{2,4}& \label{eq:feasibility_initial_problem1}\\
&\mathcal{L}(V,u_1){\rm~is~SOS},&\forall \beta_1,\beta_2\in[0.5,1]  \label{eq:feasibility_initial_problem2}
\end{eqnarray}
we have obtained a polynomial function as follows
\begin{eqnarray*}
V_0 &=&     17.6626x_1^2  -18.2644  x_1 x_2+ 16.4498 x_2^2\\
         &&          -0.1542x_1^3  +  1.7303x_1^2 x_2   -1.0845 x_1 x_2^2\\
            &&     +  0.1267x_2^3    +3.4848x_1^4   -0.8361 x_1^3 x_2\\
             &&     +4.8967 x_1^2x_2^2  + 2.3539x_2^4
\end{eqnarray*}
which, together with \eqref{eq:ex2_u1}, satisfies Assumption \ref{ass:polyall}.

Only for simulation purpose, we set $\beta_1=0.7$ and $\beta_2=0.6$. The initial condition is arbitrarily set as $x_1(0)=1$ and $x_2(0)=-2$.
The proposed online learning scheme is applied to update the control policy every four second for seven times. The exploration noise is the sum of sinusoidal waves with different frequencies, and it is turned off after the last iteration.  The suboptimal and globally stabilizing control policy is is $u_8 = [u_{8,1}, u_{8,2}]^T$ with
\begin{eqnarray*}
u_{8,1} &=&
- 0.0004x_1^3 - 0.0067x_1^2x_2 - 0.0747x_1^2\\
&& + 0.0111x_1x_2^2 + 0.0469x_1x_2 - 0.2613x_1\\
&& - 0.0377x_2^3 - 0.0575x_2^2 - 2.698x_2\\
u_{8,2} &=&
   - 0.0005x_1^3 - 0.009x_1^2x_2 - 0.101x_1^2\\
 &&  + 0.0052x_1x_2^2 - 0.1197x_1x_2 - 1.346x_1\\
 && - 0.0396x_2^3 - 0.0397x_2^2 - 3.452x_2.
\end{eqnarray*}
The associated cost function is as follows:
\begin{eqnarray*}
V_8 &=& 1.4878x_1^2  +  0.8709x_1x_2 +   4.4963x_2^2\\
&&  +0.0131 x_1^3+ 0.2491x_1^2x_2   -0.0782x_1x_2^2\\
&& +   0.0639x_2^3 +0.0012x_1^3x_2+ 0.0111x_1^2x_2^2\\
&&   -0.0123x_1x_2^3 +   0.0314 x_2^4.
\end{eqnarray*}

\begin{figure}[!t]
\centering
\includegraphics[width=0.5\textwidth]{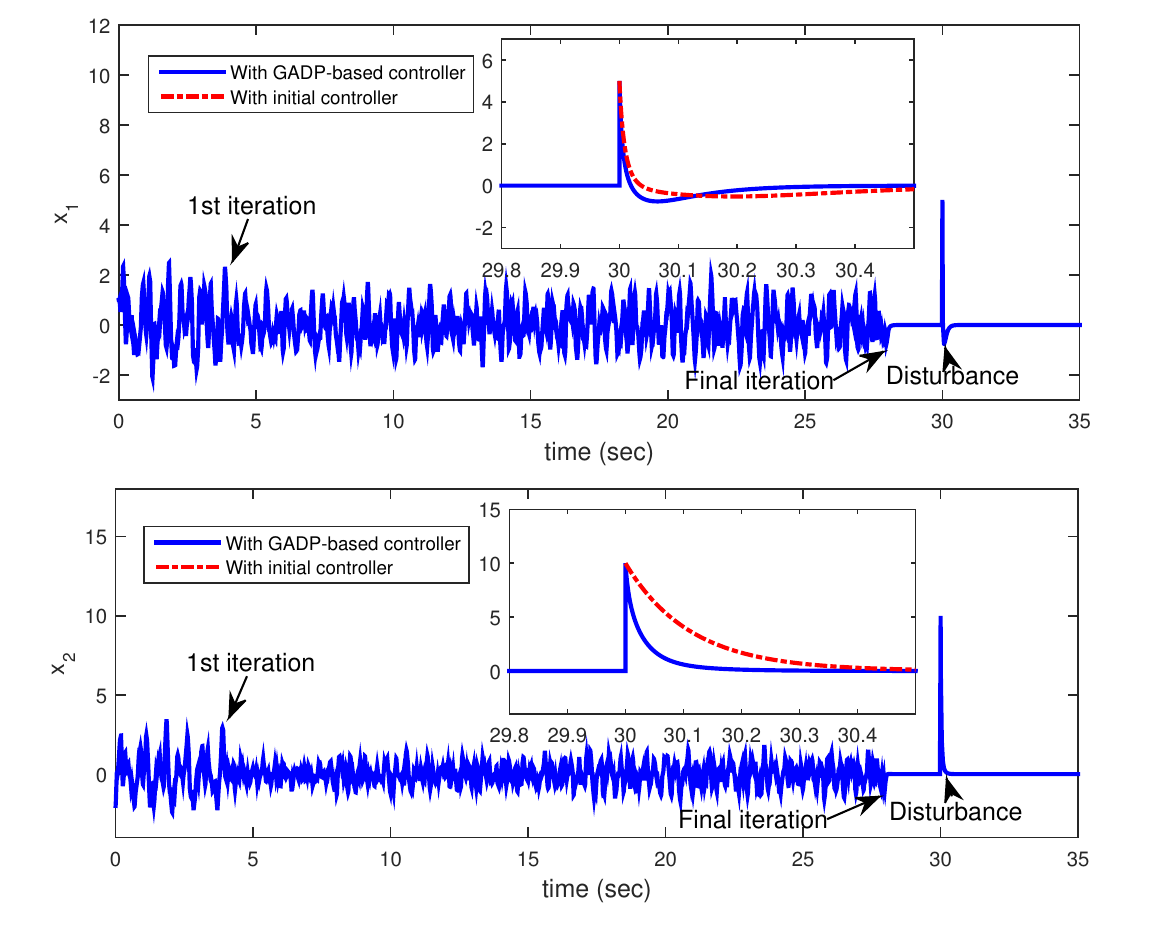}
\caption{System trajectories.}
\label{fig:ex2_state}
\end{figure}

\begin{figure}[!t]
\centering
\includegraphics[width=0.5\textwidth]{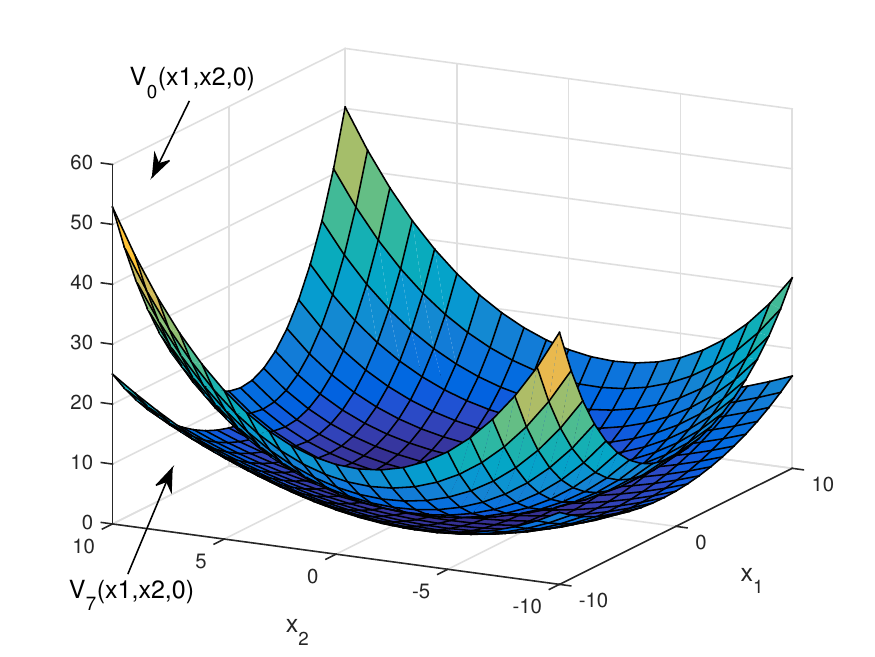}
\caption{Comparison of the cost functions.}
\label{fig:ex2_cost}
\end{figure}

In Figure \ref{fig:ex2_state}, we show the system state trajectories during the learning phase and the post-learning phase. At $t=30$, we inject an impulse disturbance through the input channel to deviate the state from the origin. Then, we compare the system response under the proposed control policy and the initial control policy given in \cite{xu2009simultaneous}. The suboptimal cost function and the original cost function are compared in Figure \ref{fig:ex2_cost}.


\subsection{An active suspension system}
Consider the quarter-car suspension model as described by the following set of differential equations \cite{gaspar2003active}.
\setlength{\arraycolsep}{5pt}
\begin{eqnarray}
\dot x_1 &=& x_2 \label{eq:suspen1}\\
\dot x_2 &=& -\frac{k_s(x_1-x_3)+k_{n}(x_1-x_3)^3}{m_b}\nonumber\\
&&-\frac{b_s(x_2-x_4) - u}{m_b}\\
\dot x_3 &=& x_4\\
\dot x_4 &=&  \frac{k_s(x_1-x_3)+k_n(x_1-x_3)^3}{m_w} \nonumber\\
&&  +\frac{b_s(x_2-x_4)+k_tx_3- u}{m_w}\label{eq:suspenn}
\end{eqnarray}
\setlength{\arraycolsep}{5pt}
where $x_1$, $x_2$, and $m_b$ denote respectively the position, velocity, and mass of the car body;  $x_3$, $x_4$, and $m_w$ represent repectively the position, velocity, and mass of the  wheel assembly; $k_t$, $k_s$, $k_n$, and $b_s$ are the tyre stiffness, the linear suspension stiffness, the nonlinear suspension stiffness, and the damping rate of the suspension.

It is assumed that $m_b\in[250,350]$, $m_w\in[55,65]$, $b_s\in[900,1100]$, $k_s\in[15000,17000]$, $k_n =k_s/10$, and $k_t\in [180000,200000]$. Then, it is easy to see that, without any control input, the system is globally asymptotically stable at the origin. Here we would like to use the proposed online learning algorithm to design an active suspension control system which reduces the following performance index
\begin{eqnarray}
J(x_0, u) = \int_0^\infty (\sum_{i=1}^4x_i^2+u^2) dt
\end{eqnarray}
and at the same time maintain global asymptotic stability. In particular, we would like to improve the system performance in the set $\Omega=\{x|x\in\mathbb{R}^4~{\rm and}~|x_1|\le 0.05, |x_2|\le 10, |x_3|\le 0.05, |x_4|\le 10, \}$.

To begin with, we first use SOSTOOLS \cite{sostools, prajna2002introducing} to obtain an initial cost function $V_0\in\mathbb{R}[x]_{2,4}$ for the polynomial system \eqref{eq:suspen1}-\eqref{eq:suspenn} with uncertain parameters, of which the range is given above. This is similar to the SOS feasibility problem described in \eqref{eq:feasibility_initial_problem1}-\eqref{eq:feasibility_initial_problem2}.

Then, we apply the proposed online learning method with $u_1 = 0$. The initial condition is arbitrarily selected. From $t=0$ to $t=120$, we apply bounded exploration noise as inputs for learning purpose, until convergence is attained after $10$ iterations. 

The suboptimal and global stabilizing control policy we obtain is
\setlength{\arraycolsep}{0.0em}
\begin{eqnarray*}
u_{10} &=&
- 3.53 x_1^3 - 1.95 x_1^2 x_2 + 10.1 x_1^2 x_3 + 1.11 x_1^2 x_4\\
&&- 4.61 \times 10^8 x_1^2 - 0.416 x_1 x_2^2 + 3.82 x_1 x_2 x_3\\
 &&+ 0.483 x_1 x_2 x_4 + 4.43 \times 10^8 x_1 x_2 - 10.4 x_1 x_3^2\\
  &&- 2.55 x_1 x_3 x_4 - 2.91 \times 10^8 x_1 x_3 - 0.174 x_1 x_4^2\\
   &&- 2.81 \times 10^8 x_1 x_4 - 49.2 x_1 - 0.0325 x_2^3 + 0.446 x_2^2 x_3\\
    &&+ 0.06 x_2^2 x_4 + 1.16 \times 10^8 x_2^2 - 1.9 x_2 x_3^2 - 0.533 x_2 x_3 x_4\\
     &&- 6.74 \times 10^8 x_2 x_3 - 0.0436 x_2 x_4^2 - 2.17 \times 10^8 x_2 x_4\\
      &&- 17.4 x_2 + 3.96 x_3^3 + 1.5 x_3^2 x_4 + 7.74 \times 10^8 x_3^2\\
       &&+ 0.241 x_3 x_4^2 + 2.62 \times 10^8 x_3 x_4 + 146.0 x_3 + 0.0135 x_4^3\\
       && + 1.16 \times 10^8 x_4^2 + 12.5 x_4
\end{eqnarray*}

\begin{figure}[!t]
\centering
\includegraphics[width=0.5\textwidth]{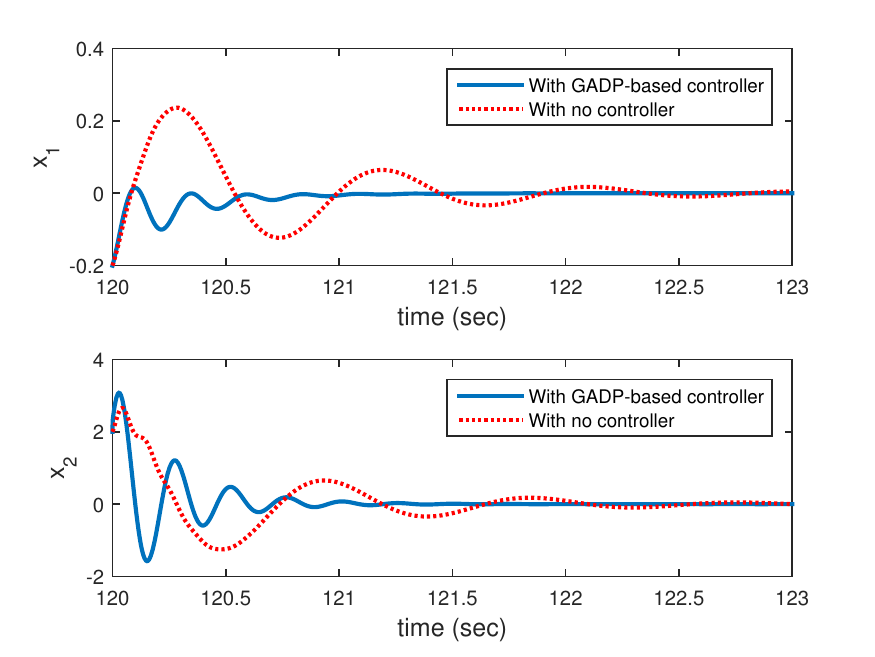}
\includegraphics[width=0.5\textwidth]{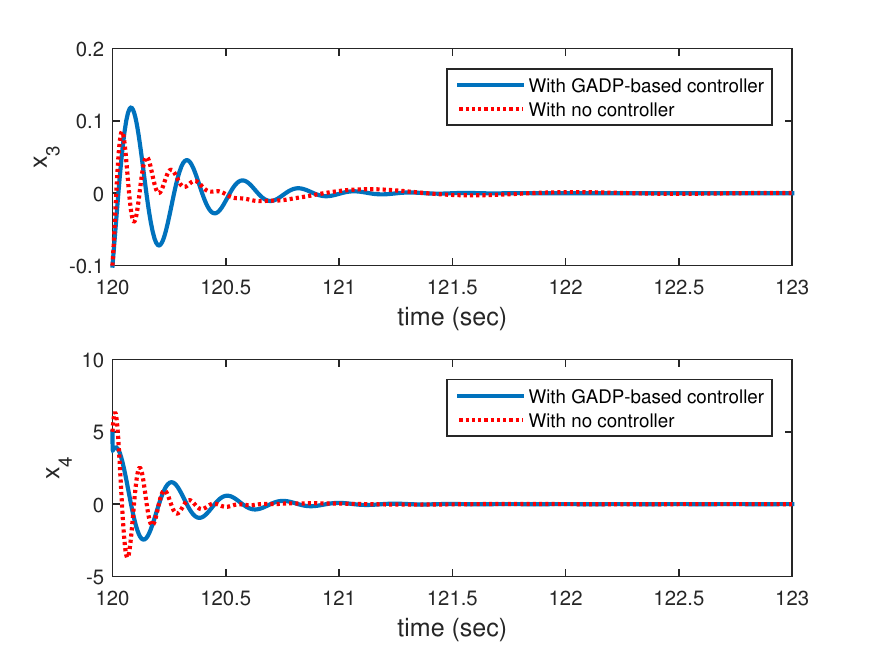}
\caption{Post-learning performance.}
\label{fig:ex3_state}
\end{figure}

At $t=120$, an impulse disturbance is simulated such that the state is deviated from the origin. Then, we compare the post-learning performance of the closed-loop system using the suboptimal control policy with the performance of the original system with no control input (see Figure \ref{fig:ex3_state}).

\begin{figure}[!t]
\centering
\includegraphics[width=0.45\textwidth]{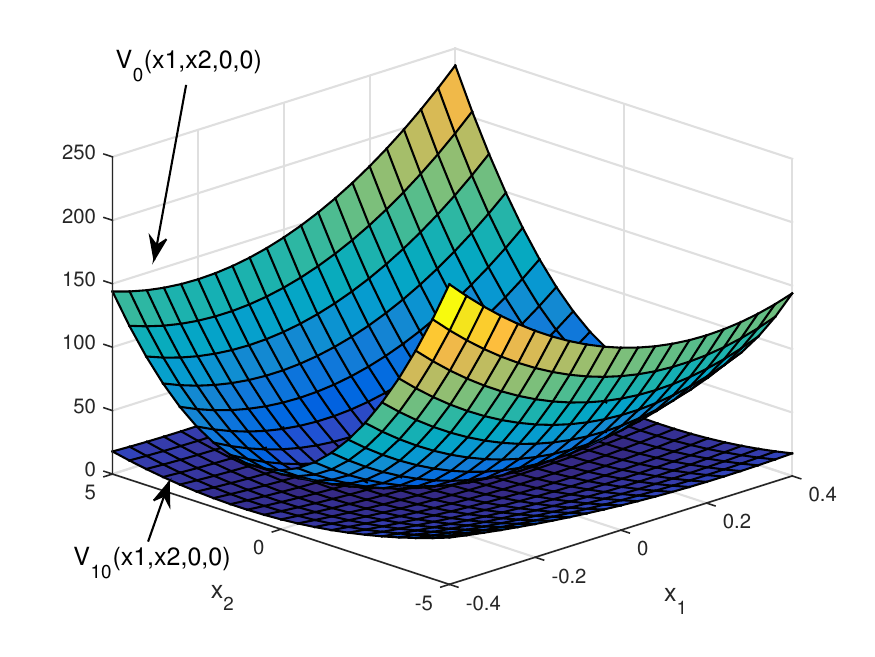}
\caption{Comparison of the cost functions.}
\label{fig:ex3_cost}
\end{figure}

To save space, we do not show their explicit forms of $V_0$ and $V_{10}$, since each of them is the summation of 65 different monomials. In Figure \ref{fig:ex3_cost}, we compare these two cost functions by restricting $x_3\equiv x_4 \equiv 0$. It can be seen that $V_{10}$ has been significantly reduced from $V_0$.
\setlength{\arraycolsep}{5pt}

\section{Conclusions}
This paper has proposed, for the first time,  a global ADP method for the data-driven (adaptive) optimal control of  nonlinear polynomial systems. In particular, a new policy iteration scheme has been developed. Different from conventional policy iteration, the new iterative technique does not attempt to solve a partial differential equation but a convex optimization problem at each iteration step. It has been shown that, this method can find a suboptimal solution to continuous-time nonlinear optimal control problems \cite{lewis2012optimal}. In addition, the resultant control policy is globally stabilizing. Also, the method can be viewed as a computational strategy to solve directly Hamilton-Jacobi inequalities, which are used in $H_{\infty}$ control problems \cite{helton1999extending, schaft1999l2}.

When the system parameters are unknown, conventional ADP methods utilize neural networks to approximate online the optimal solution, and a large number of basis functions are required to assure high approximation accuracy on some compact sets. Thus, neural-network-based ADP schemes may result in slow convergence and loss of global asymptotic stability for the closed-loop system. Here, the proposed global ADP method has overcome the two above-mentioned shortcomings, and it yields computational benefits.

It is under our current investigation to extend the proposed methodology for more general (deterministic or stochastic) nonlinear systems \cite{tao2014auto}, \cite{Horowitz_ACC2014}, and \cite{Horowitz_Preprint2014}, as well as systems with parametric and dynamic uncertainties \cite{jiang1998design, jiang2012tnnls, jiang2013robust, tao2014tie}.

\section*{Acknowledgement}
The first author would like to thank Dr. Yebin Wang, De Meng, Niao He, and Zhiyuan Weng for many helpful discussions on semidefinite programming, in the summer of 2013.


%
%

\appendices
\section{Sum-of-squares (SOS) program}

%
%
%
%
%

An SOS program is a convex optimization problem of the following form
\begin{problem}[SOS programming~\cite{blekherman2013semidefinite}]
\begin{eqnarray}
\min_y&& b^Ty\\
{\rm s.t.}&& p_i(x;y){\rm~are~SOS},~i=1,2,\cdots,k_0 \label{eq:sos_constraint}
\end{eqnarray}
where $p_i(x;y)=a_{i0}(x)+\sum_{j=1}^{n_0}a_{ij}(x)y_j$, and $a_{ij}(x)$ are given polynomials in $\mathbb{R}[x]_{0,2d}$.
\end{problem}

In \cite[p.74]{blekherman2013semidefinite}, it has been pointed out that SOS programs are in fact equivalent to semidefinite programs (SDP)  \cite{vandenberghe1996semidefinite}, \cite{blekherman2013semidefinite}.
SDP  is a broad generalization of linear programming. It concerns with the optimization of a linear function subject to linear matrix inequality constraints, and is of great theoretic and  practical interest \cite{blekherman2013semidefinite}.
The conversion from an SOS to an SDP  can be performed either manually, or automatically using, for example, the MATLAB toolbox SOSTOOLS \cite{sostools, prajna2002introducing}, YALMIP \cite{lofberg2004yalmip}, and Gloptipoly \cite{henrion2003gloptipoly}.

\section{Proof of Theorem \ref{thm:conv}}
\label{sec:proofConv}
Before proving Theorem \ref{thm:conv}, let us first give the following lemma.
\begin{lemma}
\label{le:pi}
Consider the conventional policy iteration algorithm described in \eqref{eq:GHJB} and \eqref{eq:pimp}. Suppose $u_i(x)$ is a globally stabilizing control policy and there exists $V_i(x)\in\mathcal{C}^1$ with $V(0)=0$,  such that \eqref{eq:GHJB} holds. Let $u_{i+1}$ be defined as in \eqref{eq:pimp}. Then, Under Assumption \ref{ass:HJB}, the followings are true.
\begin{enumerate}
\item $V_i(x) \ge V^{\rm o}(x)$;
\item for any $V_{i-1}\in\mathcal{P}$, such that $\mathcal{L}(V_{i-1},u_i)\ge 0$, we have $V_{i}\le V_{i-1}$;
\item $\mathcal{L}(V_{i},u_{i+1})\ge 0$.
\end{enumerate}
\end{lemma}

\begin{IEEEproof}
1) Under Assumption \ref{ass:HJB}, we have
\begin{eqnarray*}
0&=&\mathcal{L}(V^{\rm o}, u^{\rm o})-\mathcal{L}(V_i, u_i)\\
&=&(\nabla  V_i-\nabla V^{\rm o})^T(f+gu_i)+r(x,{u_i})\\
&&- (\nabla V^{\rm o})^Tg(u^{\rm o}-u_i)-r(x,u^{\rm o})\\
&=&(\nabla V_i-\nabla V^{\rm o})^T(f+gu_i)+|{u_i}-u^{\rm o}|^2_R
\end{eqnarray*}

Therefore, for any $x_0\in\mathbb{R}^n$, along the trajectories of system \eqref{eq:sys1} with $u=u_i$ and $x(0)=x_0$, we have
\begin{eqnarray}
V_i(x_0)-V^{\rm o}(x_0) &=& \int_0^{T}|{u_i}-u^{\rm o}|^2_Rdt\nonumber \\
& & + V_i(x(T)) - V^{\rm o}(x(T)) \label{eq:ViVo}
\end{eqnarray}

Since $u_i$ is globally stabilizing, we know $\lim\limits_{T\rightarrow+\infty}V_i(x(T))=0$ and $\lim\limits_{T\rightarrow+\infty}V^{\rm o}(x(T))=0$. Hence, letting $T\rightarrow +\infty$, from \eqref{eq:ViVo} it follows that $V_i(x_0)\ge V^{\rm o}(x_0)$. Since $x_0$ is arbitrarily selected, we have  $V_i(x)\ge V^{\rm o}(x)$, $\forall x\in\mathbb{R}^n$.

2) Let $q_i(x)$ be a positive semidefinite function, such that
\begin{eqnarray}
\mathcal{L}(V_{i-1}(x),u_i(x))=q_i(x),~~\forall x\in \mathbb{R}^n.
\end{eqnarray}

Therefore,
\begin{eqnarray}
(\nabla V_{i-1}-\nabla V_i)^T(f+gu_i)+q_i(x) = 0.
\end{eqnarray}

Similar as in 1),  along the trajectories of system \eqref{eq:sys1} with $u=u_i$ and $x(0)=x_0$, we can show
\begin{eqnarray}
V_{i-1}(x_0)-V_i(x_0) &=& \int_0^{\infty}q_i(x)dt
\end{eqnarray}

Hence, $V_{i-1}(x)\le V_{i}(x)$, $\forall x\in\mathbb{R}^n$.

3) By definition,
\begin{eqnarray}
&&\mathcal{L}(V_{i},u_{i+1}) \nonumber\\
&=& -\nabla V_{i}^T\left(f+gu_{i+1}\right)-q-|u_{i+1}|_R^2\nonumber\\
&=& -\nabla V_{i}^T\left(f+gu_{i}\right)-q-|u_{i}|_R^2\nonumber\\
&&  -\nabla V_{i}^Tg\left(u_{i+1}-u_{i}\right)-|u_{i+1}|_R^2+|u_{i}|_R^2\nonumber\\
&=& 2u_{i+1}^TR\left(u_{i+1}-u_{i}\right)-|u_{i+1}|_R^2+|u_{i}|_R^2\nonumber\\
&=& -2u_{i+1}^TRu_{i}+|u_{i+1}|_R^2+|u_{i}|_R^2\nonumber\\
&\ge& 0
\end{eqnarray}

The proof is complete.\\
\end{IEEEproof}

\noindent {\it Proof of Theorem \ref{thm:conv}}

We first prove 1) and 2) by induction. To be more specific, we will show that 1) and 2) are true and $V_i\in\mathcal{P}$, for all $i=0,1,\cdots$.

i) If $i=1$, by Assumption \ref{ass:controlability} and Lemma \ref{le:pi} 1), we immediately know 1) and 2) hold. In addition, by Assumptions \ref{ass:controlability} and \ref{ass:HJB}, we have $V^{\rm o}\in\mathcal{P}$ and $V_0\in\mathcal{P}$. Therefore, $V_i\in\mathcal{P}$.

ii) Suppose 1) and 2) hold for $i=j>1$, and $V_j\in\mathcal{P}$. We show  1) and 2) also hold for $i=j+1$, and $V_{j+1}\in\mathcal{P}$.

Indeed, since $V^{\rm o}\in\mathcal{P}$ and $V_j\in\mathcal{P}$. By the induction assumption, we know $V_{j+1}\in\mathcal{P}$.

Next, by Lemma \ref{le:pi} 3), we have $\mathcal{L}(V_{j+1},u_{j+2})\ge0$. As a result, along the solutions of system \eqref{eq:sys1} with $u=u_{j+2}$, we have
\begin{eqnarray}
\dot V_{j+1}(x)\le -q(x).
\end{eqnarray}
Notice that, since $V_{j+1}\in\mathcal{P}$, it is a well-defined Lyapunov function for the closed-loop system \eqref{eq:sys1} with $u=u_{j+2}$.  Therefore, $u_{j+2}$ is globally stabilizing, i.e., 2) holds for $i=j+1$.

Then, by Lemma \ref{le:pi} 2), we have $V_{j+2}\le V_{j+1}$. Together with the induction Assumption, it follows that
$V^{\rm o}\le V_{j+2}\le V_{j+1}$. Hence, 1) holds for $i=j+1$.
\vspace{2mm}

Now, let us prove 3). If such a pair $(V^*,u^*)$ exists, we immediately know $u^*=-\frac{1}{2}R^{-1}g^T\nabla V^*$. Hence,
\begin{eqnarray}
\mathcal{H}(V^*)=\mathcal{L}(V^*,u^*)=0.
\end{eqnarray}
Also, since $V^{\rm o}\le V^*\le V_0$, $V^*\in\mathcal{P}$.
However, as discussed in Section \ref{sec:optnstb}, solution to the HJB equation \eqref{eq:HJB} must be unique. Hence, $V^*=V^{\rm o}$ and $u^*=u^{\rm o}$.

The proof is complete.
\QED
\bibliographystyle{IEEEtranS}
\bibliography{IEEEabrv,gadpref}                                        
\end{document}